\documentclass[a4paper,12pt,leqno]{article}
\usepackage{latexsym}
\usepackage[all]{xy}
\usepackage{mathrsfs} 
\usepackage{amssymb} 
\usepackage{amsmath} 
\usepackage{theorem}

\numberwithin{equation}{section}

\newcommand{\owari}{\hfill$\square$}

\theoremstyle{break}
\newtheorem{theorem}{Theorem}[section]

\newtheorem{cor}[theorem]{Corollary}
\newtheorem{lemma}[theorem]{Lemma}

\newtheorem{rem}[theorem]{Remark}

\newtheorem{examples}[theorem]{Example}

\title{Combinatorial data of a free arrangement and the Terao conjecture }
\author{
Tran Quoc Cong
}

\begin{document}

\maketitle

\begin{abstract}
Let $\mathcal{A}$ be an hyperplane arrangement in a real or complex vector space. We denote by $D(\mathcal{A})$ the module of logarithmic derivations. In this paper, we present a combinatorial structure of generators of $D(\mathcal{A}).$ This structure permits us to detect the relationship between the combinatorial determined property and the singularity of vector field. Consequently, by using only combinatorial data, we have a basis of the module in free case and that yields a proof  for the Terao's conjecture. We also verify the example of Ziegler and give a sufficient condition on combinatorial determined property of generators.
\end{abstract}
 
\section{Introduction}
We consider an arrangement $\mathcal{A}$ consisting of a finite collection of linear hyperplanes in a real or complex vector space $V$. $\mathcal{A}$ is {\it free} if its module of logarithmic derivations is a free $\mathrm{S}\mathrm{y}\mathrm{m}(V^{*})$-module. In free arrangement context, there is a central problem is that called Terao conjecture, which asserts the freeness of $\mathcal{A}$ depends only on the combinatorics of the arrangement. And more general problem, to determine whether a given arrangement is free or not, is also a very interesting problem. These problems served as a major motivation for many results in arrangements of hyperplanes. Their advances are contributed by many authors (see \cite{Sa1}, \cite{T0}, \cite{T3}, \cite{Yo}, \cite{Yu2}, \cite{Zi}, etc).
\\
In order to describe module $D(\mathcal{A})$ of any arrangement in general, one has an unique tool up to now: Gr\"{o}bner basis (see \cite{BC}). It is easy to use this tool for computer but there are many computations which would be extremely intractable to do by hand. Therefore, it is very difficult to control combinatorial data as well as freeness. In this paper, we introduce a new approach to the problems without using Gr\"{o}bner basis: to investigate its system of equations instead of module $D(\mathcal{A})$.\\
The organization of this is as follows. In Section 2 we give some results that they will lead us to a system of equations describing combinatorial data of $D(\mathcal{A})$; to recognize the essence of combinatorial determined property lies entirely in the non-homogeneous parts. In Section 3 we present a combinatorial structure of generators, verify the example of Ziegler and give a sufficient condition on combinatorial determined property as a generalization. Consequently, when it is free, derivations $\theta_i$ of a basis of $D(\mathcal{A})$ satisfy that sufficient condition. They yields a proof  for the Terao's conjecture in the last corollary. 
\medskip

\noindent
\textbf{Acknowledgements}. 
The author is deeply grateful to Nguyen Viet Dung for his advices and discussions.
He wishes to express his sincere thanks to Sergey Yuzvinsky for his comments on the example of Ziegler.
\section{Combinatorial data of module $D(\mathcal{A})$}
Firstly, we review some basic concepts concerning free arrangements. Our standard reference is \cite{OT}. \\
Let $\mathcal{A}$  be an hyperplane arrangement in $\ell-$dimensional vector space $V$ over $\mathbb{K}$. $L(\mathcal{A})$ is the lattice and $M(\mathcal{A})$ is the complement of $\mathcal{A}$. In order to investigate the freeness, it is sufficient to suppose that $\mathcal{A}$ is essential (see \cite{Zi}). For each hyperplane $H_i\in\mathcal{A}$, choose a defining equation $\alpha_i\in V^*$ and put $Q(\mathcal{A})=\prod_{1\leq i\leq n}\alpha_i$. We identify $S:=\mathrm{S}\mathrm{y}\mathrm{m}(V^{*})$ with the polynomial algebra $K[x_1,\cdots, x_\ell ]$ and denote by $Der_K$ be the module of derivations of $S$ over $K.$ The arrangement $\mathcal{A}$ is \textit{free} if  its module of logarithmic derivations 
$$ D(\mathcal{A})=\{\theta\in Der_K| \theta Q\in QS\}$$ is a free $S-$module. To grasp the combinatorial essence of every $\theta=\sum_{i=1}^\ell p_i.D_i\in D(\mathcal{A})$, we only need to work on its vector field $\bar{\theta}=(p_1, p_2,\cdots, p_\ell )$, where $p_1, p_2, ..., p_\ell\in S.$\\
Next, we give some results that they will lead us to an alternative approach for module $D(\mathcal{A})$.
\begin{theorem}
For any $\theta\in  D(\mathcal{A})$, with some polynomials $k_1, \cdots, k_n$, then $\bar{\theta }$ satisfies the following equation:
$$\begin{bmatrix}
a_{11}&a_{12}&...&a_{1\ell}\\ 
a_{21}&a_{22}&...&a_{2\ell}\\ 
...&...&...&...\\ 
a_{n1}&a_{n2}&...&a_{n\ell}\\ 
\end{bmatrix}.\begin{bmatrix}
p_1\\ 
p_2\\ 
...\\ 
p_\ell\\ 
\end{bmatrix}=\begin{bmatrix}
k_1.\alpha_1\\ 
k_2.\alpha_2\\ 
...\\ 
k_n.\alpha_n\\ 
\end{bmatrix},$$ 
where $\alpha_i= a_{i1}x_1 +\cdots + a_{i\ell }x_\ell, \forall i=\overline{1,n.}$
\end{theorem}
\textbf{Proof}. 
Since $\theta\in D(\mathcal{A})$, we may write $\theta Q = gQ$ with some polynomial $g\in S.$ Thus 
\begin{multline*} 
\theta Q=\theta \left[\prod_{i=1}^n\sum_{j=1}^\ell a_{ij}x_j\right]=
(\sum_{s=1}^\ell p_s D_s)\left[\prod_{i=1}^n\sum_{j=1}^\ell a_{ij}x_j\right]\\
=\frac{Q}{\alpha_1}\sum_{j=1}^\ell a_{1j}p_j+\cdots+\frac{Q}{\alpha_n}\sum_{j=1}^\ell a_{nj}p_j=gQ.
\end{multline*}
This shows that 
$$ \frac{1}{\alpha_1}\sum_{j=1}^\ell a_{1j}p_j+\cdots+\frac{1}{\alpha_n}\sum_{j=1}^\ell a_{nj}p_j=g.$$
Since the linear forms $\alpha_i$ are coprime, $\alpha_i$ divides $\sum_{j=1}^\ell a_{ij}p_j$, for all $1\leq i\leq n.$ This means that there exists polynomials $k_1,\cdots,k_n$ such that $\sum_{j=1}^\ell a_{ij}p_j =k_i \alpha_i, 1\leq i\leq n.$ Thus the theorem is proved. \owari

Without loss of generality, we can assume that the first $\ell$ normal vectors $n_{H_1}, \cdots, n_{H_\ell}$ of the linear forms $\alpha_1, \cdots, \alpha_\ell$ are linearly independent. Hence, the linear forms $\alpha_j$'s are expressible uniquely in terms of the linear forms $\alpha_1, \cdots, \alpha_\ell$
$$ \alpha_j=\sum_{i=1}^\ell a_{ij}\alpha_i, \ell+1\leq j\leq n.$$
\begin{cor}
For any $\theta\in D(\mathcal{A})$, there is a $\ell-$tuple $[k_1, k_2,\cdots, k_\ell]$ satisfying the following equations
\begin{equation}
  k_j\alpha_j=\sum_{i=1}^\ell a_{ij}\alpha_ik_i, \ell+1\leq j\leq n.
\end{equation}
\end{cor}
\textbf{Proof}. 
By theorem $2.1$, $\theta$ satisfies 
$$\begin{bmatrix}
a_{11}&a_{12}&...&a_{1\ell}\\ 
a_{21}&a_{22}&...&a_{2\ell}\\ 
...&...&...&...\\ 
a_{n1}&a_{n2}&...&a_{n\ell}\\ 
\end{bmatrix}.\begin{bmatrix}
p_1\\ 
p_2\\ 
...\\ 
p_\ell\\ 
\end{bmatrix}=\begin{bmatrix}
k_1.\alpha_1\\ 
k_2.\alpha_2\\ 
...\\ 
k_n.\alpha_n\\ 
\end{bmatrix}.$$
The system of equations is obtained by elementary row operations.\owari\\
The defining polynomial of arrangement $\mathcal{A}$ is called {\it canonical } if it has form $Q(\mathcal{A})=x_1\cdots x_\ell.\alpha_{\ell+1}\cdots\alpha_n$. In order to investigate the freeness of arrangements, we need only to consider their defining polynomials in canonical form.
\begin{lemma}
By changing suitable coordinates, we can obtain a defining polynomial of a given arrangement in canonical form without changing its freeness and lattice.
\end{lemma}
\textbf{Proof}. 
Let  $\mathcal{A}$ be an arrangement with defining polynomial $Q(\mathcal{A})=\alpha_1\cdots\alpha_n$ in polynomial ring $S=K[x_1,\cdots, x_\ell ]$. Without loss of generality, we can assume that the first $\ell$ linear forms $\alpha_1, \cdots, \alpha_\ell$ are linearly independent. We change our coordinate system by taking $X_i=\alpha_i, \forall i=1,\cdots, \ell$ or $$\begin{bmatrix}
X_1\\ 
X_2\\ 
...\\ 
X_\ell\\ 
\end{bmatrix}=A.\begin{bmatrix}
x_1\\ 
x_2\\ 
...\\ 
x_\ell\\ 
\end{bmatrix}.$$

The defining polynomial becomes $Q=X_1\cdots X_\ell.\beta_{\ell+1}\cdots\beta_n$ in polynomial ring $S'=K[X_1,\cdots,X_\ell ]$. Next, we consider an arrangement $\mathcal{B}$ with defining polynomial $Q(\mathcal{B})=x_1\cdots x_\ell\beta_{\ell+1}\cdots\beta_n$ in polynomial ring $S=K[x_1,\cdots,x_\ell ]$. Since the changing is non-degenerate, $\mathcal{A}$ and $\mathcal{B}$ have the same lattice.\\
It remain to prove that $\mathcal{A}$ and $\mathcal{B}$ have the same freeness.\\
If $\mathcal{A}$ is free, there exists a basis $\theta_1,\cdots,\theta_\ell$ of module $D(\mathcal{A})$ such that $$det M(\theta_1,\cdots,\theta_\ell)\big{|}_{(x_1,\cdots,x_\ell)}=c.Q(\mathcal{A}).$$
Therefore, one has 
\begin{multline*} 
det M(\theta_1,\cdots,\theta_\ell)\big{|}_{(X_1,\cdots,X_\ell)}=det A. det M(\theta_1,\cdots,\theta_\ell)\big{|}_{(x_1,\cdots,x_\ell)}
\\=c.det A.Q(\mathcal{A})=c.det A.Q.
\end{multline*} 
It means that $\mathcal{A}$ has the same freeness to some arrangement $\mathcal{C}$ having defining polynomial $Q$ in polynomial ring $S'=K[X_1,\cdots,X_\ell ]$. Hence, $\mathcal{A}$ and $\mathcal{B}$ have the same freeness.
\owari\\

From now on, we only consider arrangements with defining polynomials in canonical form. In this case, system (2.1) becomes \begin{equation}
  k_j\alpha_j=\sum_{i=1}^\ell k_i(a_{ij}x_i), \ell+1\leq j\leq n.
\end{equation}
A {\it solution} is a $n$-tuple $\theta=[k_1, k_2,\cdots,k_n]$ of polynomials satisfying $(2.2)$. And $\theta$ is defined uniquely by the first $\ell$-tuple $[k_1, k_2,\cdots,k_\ell]$. Therefore, if no confusion is possible, we still denote by $\theta=[k_1, k_2,\cdots,k_\ell]$ to be a solution.  For each $j=\overline{\ell+1,n}$, we call $\sum_{i=1}^\ell k_i(a_{ij}x_i)=0$ {\it j-th homogeneous equation} and $k_j\alpha_j=\sum_{i=1}^\ell k_i(a_{ij}x_i)$ {\it j-th non-homogeneous equation}.
\begin{lemma}
Set of solutions of j-th equation is generated by the following system of canonical generators:
\begin{center}$G_j=$
$\begin{cases}
e=[1, 1,\cdots, 1],&\\
e_r=[0,\cdots, 0, 1, 0,\cdots, 0], \ {\text if} \ \ a_{jr}=0,&\\
a_{jt}x_t.e_s -a_{js}x_s.e_t, \ {\text if}\ \ a_{jt}.a_{js}\neq 0,&
\end{cases}$
\end{center}
where $e_r$ is $r$-th identity vector.
\end{lemma}
\textbf{Proof}. 
One has two cases:\\
{\it (i) Homogeneous case.} We will construct a system of generators for the equation $\sum_{i=1}^\ell k_i(a_{ij}x_i)=0$. Indeed, the set $G_j \setminus \{e\}$ is a system of generators for {\it linear syzygies}.

\vspace{5mm}

\noindent
{\it (ii) Non-homogeneous case.} Any non-homogeneous solution $\theta=[k_1,\cdots, k_\ell]$ can rewrite as $\theta=k_j.e -\gamma,$ where $\gamma$ is a solution of homogeneous equation.\\
Therefore, $G_j$ is a system of generators of $j-$th equation.
\owari
\vspace{5mm}

\section{Combinatorial structure on generators}
For any $\theta\in D(\mathcal{A})$, its $\bar{\theta}$ can rewrite $\bar{\theta}=\sum_{k=1}^M m_k.v_k,$ where $M$ is the number of monomials $m_k$ of $\bar{\theta}$, $v_k\in\mathbb{K}^\ell.$ We will describe constraints arising in sets $\mathscr{C}_j=\{\alpha_j(v_k)| k=\overline{1, M}\},j=\overline{1,n}.$

\begin{theorem}[Combinatorial Structure]
(i) interior constraints : there exists a subset $I\subseteq [M], c_{ij}\in\mathscr{C}_j$ and coefficents $b_i\in\mathbb{K}$ such that $\sum_{i\in I}b_i. c_{ij}=0.$\\
(ii) exterior constraints : If $\sum_{j\in J}b_j.\alpha_j=0.$ then  $\sum_{j\in J}b_j. c_{ij}=0.$
\end{theorem}
\textbf{Proof}.
(ii) It is followed from $\sum_{j\in J}b_j.\alpha_j=0.$\\
(i) According to lemma $2.4$, we have $\bar{\theta}=\bar{\theta}_h^i +\bar{\theta}_{nh}^i$, where $\bar{\theta}_h^i (\bar{\theta}_{nh}^i)$ is homogeneous(non-homogeneous) vector field with respect to $H_i.$ Moreover, $\bar{\theta}_{nh}^i=\sum_k m'_k(x_1,\cdots,x_\ell).$ Therefore, \\
$$ \bar{\theta}=\bar{\theta}_h^i +\sum_k m'_k(x_1,\cdots,x_\ell). $$
It yields $$ \alpha_i(\bar{\theta})=\alpha_i(\bar{\theta}_h^i )+\alpha_i(\sum_k m'_k(x_1,\cdots,x_\ell)). $$
$$ =0+\sum_k m'_k\alpha_i(x_1,\cdots,x_\ell). $$
We have $$ \sum_{k=1}^M m_k.\alpha_i(v_k)=\sum_k m'_k\alpha_i.$$
Substitute $x=(x_1,\cdots,x_\ell)\in H_i$ into both sides and assume that monomials $m_k (k\notin I)$ vanish. Hence
$ \sum_{i\in I}b_i. c_{ij}=0. $ \owari\\ 
In addition, there are mysterious relations among $c_{ij} ( (i,j)\in [M]\times [n])$ as in the following examples. We will call them the {\it hidden contraints}. Assume that $\alpha_{i_1},\cdots\alpha_{i_\ell}$ are linearly independent. We will call the vector field $(q_{i_1},\cdots,q_{i_\ell})$ {\it associated vector field of $\theta$}, where $q_{i_k}:=\frac{\theta(\alpha_{i_k})}{\alpha_{i_k}}, (i_k=\overline{i_1,i_\ell}).$ If $(q_{i_1},\cdots,q_{i_\ell})$ has  a critical point $c\in M(\mathcal{A})$, there will be a {\it hidden constraint:}
$$ \sum_{t=i_1}^{i_\ell}\sum_{k=1}^M m_k\big|_{x=c}.\alpha_t(v_k)=0.$$  
\begin{examples}
Hidden constraint occurs in the example $8.7$ of \cite{Zi}, the remark $3.9$ of \cite{ADS} and the example $3.4.3$ of \cite{DS}.
\end{examples}
(i)Ziegler's example :\\
In \cite{Zi}, the author proved that two arrangements $X_1, X_2 (L(X_1)\cong L(X_2))$ have different degree sequences. In particular, $D(X_2)$ has a derivation $\theta_z $ of degree $5$ but $D(X_1)$ has no derivation of degree $5$ except  to $p\theta_E$.
Let consider arrangement $X_2$ with defining polynomial : {\footnotesize $Q(X_2)=xyz\alpha_4\alpha_5\alpha_6\alpha_7\alpha_8\alpha_9 =$} \\
{\footnotesize $=xyz(x+y-z)(x-y+z)(2x-2y+z)(2x-y-2z)(2x+y+z)(2x-y-z).$}\\
\\
$ \theta_z${\footnotesize $=x(4x^3y+2x^2y^2-4xy^3-2y^4+36x^3z-84x^2yz+8xy^2z+28y^3z+18x^2z^2-24xyz^2+ 
44y^2z^2-36xz^3-4yz^3-18z^4)D_x \ \  \ \ + \ \  \ \ y(8x^4+4x^3y-8x^2y^2-4xy^3-40x^3z+6x^2yz+6xy^2z+
16y^3z-54x^2z^2+84xyz^2+8y^2z^2-38xz^3-40yz^3+16z^4)D_y \ \  \ \ + \ \  \ \  z(72x^4-128x^3y+10x^2y^2+
50xy^3-16y^4+36x^3z+6x^2yz+4xy^2z-8y^3z-72x^2z^2+30xyz^2+40y^2z^2-36xz^3-16yz^3)D_z. $}\\
\\
These arrangements are different from free arrangements at crucial point: their associated vector fields $(q_1,q_2,q_3)$ have singularities in the complement.\\
\\
$q_1:=\frac{\theta_z(\alpha_4)}{\alpha_4}=${\footnotesize $12x^3y-6x^2y^2-6xy^3-36x^3z+52x^2yz-54xy^2z+32y^3z-54x^2z^2+22xyz^2+48y^2z^2-18xz^3-32yz^3$},\\
 $q_2:=\frac{\theta_z(\alpha_5)}{\alpha_5}=${\footnotesize $-4x^3y-6x^2y^2-2xy^3+108x^3z-60x^2yz-42xy^2z+32y^3z-54x^2z^2+42xyz^2+48y^2z^2-54xz^3-32yz^3$},\\
 $q_3:=\frac{\theta_z(\alpha_6)}{\alpha_6}=${\footnotesize $-4x^3y-6x^2y^2-2xy^3+72x^3z-34x^2yz-24xy^2z+24y^3z+50xyz^2+24y^2z^2-72xz^3-48yz^3$}.\\
\\
$(q_1, q_2, q_3)$ has a critical point $c=(2, 3, -1)\in M(X_2).$\\
This critical point produces a {\it hidden constraint} among $c_{ij}, ((i,j)\in [M]\times [n] )$ as follows
$$ 48c_{4,1}+\cdots + 3c_{4,18}\cdots + 48c_{6,1}+\cdots + 3c_{6,18}=0.$$
(ii) Abe-Dimca's example :\\
\\
{\footnotesize $Q(\mathcal{A}_1)=xy(x-y-z)(x-y+z)(2x+y-2z)(x+3y-3z)(3x+2y+3z)(x+5y+5z)(7x-4y-z),$}\\
\\
{\footnotesize $Q(\mathcal{A}_2)=xy(4x-5y-5z)(x-y+z)(16x+13y-20z)(x+3y-3z)(3x+2y+3z)(x+5y+5z)(7x-4y-z).$}\\
\\
These arrangements have different degree sequences. It is easy to find an hidden constraint by the same way as in (i).\\
\\
(iii) Denham-Steiner's example :\\
{\footnotesize $$Q(\mathcal{A})=x(x+y)(x+z)(x+t)(x+z+t),$$}
{\footnotesize $$Q(\mathcal{A'})=x(x+y)(x+y+z)(x+t)(x+y+z+t).$$}
It is simple to verify that $\mathcal{A'}$ is obtained from $\mathcal{A}$ by using a matrix $P\in GL(\mathbb{K}^4)$ without fixing $V(x,x+z,x+t,x+z+t).$ These arrangements have the same degree sequences $\{1,1,2,2,2\}$, but $D(\mathcal{A})$ is not isomorphic to $D(\mathcal{A'})$. This indicates that the degree sequence is insufficient to fully characterize the module $D(\mathcal{A})$. There is a natural question that now arises : what is the difference between $D(\mathcal{A})$ and $D(\mathcal{A'})$ in this case ?\\ 
To answer this question, we need to find an homorphism $f : D(\mathcal{A})\rightarrow D(\mathcal{A'})$ that is both closest to the lattices isomorphism and a map that preserves the degree sequence relying on given generators. Since $f$ preserves the degree sequence, $f$ is surjective and $Ker(f)\neq 0.$ The difference here is $Ker(f).$   
$$Ker(f)=\{S.\theta|\theta=-(x^2y+xyt)D_x + (x^2y+xyt)D_y + (x^2y+xyt)D_z\}.$$
We can use its zero $(1,0,1,1)\in M(\mathcal{A})$ to find an hidden constraint. Therefore, hidden constraints occur in this difference.\\ 
 Suppose that $\mathcal{A, A'}$ are two arrangements with $L(\mathcal{A})\cong L(\mathcal{A'})$. $\mathscr{C}_\theta(\mathcal{A})$ is the vector space generated by all of interior constraints and exterior constraints of $\theta$ (by considering $c_{ij}$ as indeterminates and $m_k$ as coefficents). $\mathscr{C}(\mathcal{A}')$ is the vector space received from $\mathscr{C}_\theta(\mathcal{A})$ by replacing $\alpha_i$ with $\alpha_i'.$
\begin{theorem}[Sufficient Condition]
If $\theta$ has no hidden constraints, then $\theta$ is combinatorially determined.
\end{theorem}
\textbf{Proof}. 
We will prove that $\theta$ is combinatorially determined by showing that there exist a derivation $\theta'\in D(\mathcal{A'})$ such that $\mathscr{C}_\theta(\mathcal{A})\cong \mathscr{C}(\mathcal{A}')\cong\mathscr{C}_\theta'(\mathcal{A'}).$\\
Since $\theta\in D(\mathcal{A})$, \begin{equation}
\bar{\theta}=\sum_{t=1}^M m_t.v_t=\bar{\theta}_h^i +\bar{\theta}_{nh}^i, i=\overline{\ell+1,n}.\end{equation}
Denote by $\bar{\theta}_v=(v_1, v_2,\cdots, v_M)^t$ the column vector with $M\times\ell$ entries.\\
By using lemma $2.4$, one has $m_t.v_t=m_t\sum_s r_{st}.g_{st},r_{st}\in\mathbb{K}, g_{st}\in\mathbb{K}^{\ell}.$ Note that if $m_t$ and $m_u$ occur in the same $m.\bar{\theta_E},$ then $r_{st} = r_{s'u}$ for some $s, s'.$\\
Therefore, we can rewrite $(3.1)$ : $$\bar{\theta}_v= A_j.(r_{11},\cdots,r_{k_11},\cdots,r_{k_1M}\cdots,r_{k_MM})^t, j=\overline{\ell+1,n}$$
It means that $\bar{\theta}_v$ is related to a homogeneous system of linear equations : $$A_{\ell+1}.X_{\ell+1}=\cdots =A_n.X_n.$$
Consider a family of arrangements with defining polynomials : $$Q(\mathcal{A}_j)=x_1.x_2\cdots x_\ell.\alpha_j ,j=\overline{\ell+1,n}.$$
There is a natural isomorphism $h_j$ between $H_j\in\mathcal{A}$ and $H'_j\in\mathcal{A}'$ as follows:
$$ (x_1,\cdots,x_\ell)\mapsto (\frac{a_1}{a_1'}x_1,\cdots,\frac{a_\ell}{a_\ell'}x_\ell) $$
(considering $\frac{a_i}{a_i'}=1$ if $a_i=a_i'=0$), $j=\overline{\ell+1,n}$. Then it induce a module isomorphism $h_j^{\star} : D(\mathcal{A}_j)\rightarrow D(\mathcal{A'}_j).$ Based on these module isomorphisms, we can construct the same homogeneous system of linear equations :
$$A'_{\ell+1}.X_{\ell+1}=\cdots =A'_n.X_n,$$
such that $rank(A_j)=rank(A'_j), \forall j=\overline{\ell+1,n}.$ Since $\theta$ has no hidden constraints, solution spaces of homogeneous systems are isomorphism. Therefore, there exist a $\theta'\in D(\mathcal{A'})$ having the same combinatorial structure.
\owari\\
Note that the isomorphism $h_j : H_j\rightarrow H'_j$ is fixing the set $(\bigcap_{i\in I}{x_i})\cap H_j.$ Thus, it differs from the example $3.4.3$ of \cite{DS}. If we replace $\mathbb{K}$ with $\mathbb{K}'$, then the space $\mathscr{C}(\mathcal{A}')$ can be deformed. Therefore, there will be phenomena as in \cite{Z0}. Moreover, one does not consider an isomorphism between $\mathscr{C}_\theta(\mathcal{A})$ and $\mathscr{C}(\mathcal{A}')$ in a situation like that.
\begin{rem}
Assume that $\mathcal{A}$ is free and $\{\theta_E,\theta_2,\cdots,\theta_i,\cdots,\theta_\ell \}$ is a basis of $D(\mathcal{A}).$ \\
(i) $\{\theta_E,\theta_2,\cdots,\theta_i+p\theta_E,\cdots,\theta_\ell \}$ is also a basis, for any polynomial $p.$\\
(ii) The associated vector fields of $\theta_i$ have no critical points in the complement.
\end{rem}
\textbf{Proof}. 
(i) The proof is a routine computation.\\
(ii) We have  $$det M(\theta_E, \theta_2,\cdots,\theta_\ell) =cQ(\mathcal{A}).$$
Using proposition $4.12$ in \cite{OT} (see page $103$)
$$\begin{vmatrix}
..&...&\theta_j(\alpha_{i_1})&...&...\\ 
..&...&\theta_j(\alpha_{i_2})&...&...\\ 
..&...&...&...&...\\ 
..&...&\theta_j(\alpha_{i_\ell})&...&...\\ 
\end{vmatrix}=cQ(\mathcal{A}).$$
Therefore 
$$\begin{vmatrix}
..&...&\alpha_{i_1}q_1&...&...\\ 
..&...&\alpha_{i_2}q_2&...&...\\ 
..&...&...&...&...\\ 
..&...&\alpha_{i_\ell}q_\ell&...&...\\ 
\end{vmatrix}=cQ(\mathcal{A}).$$
It yields $$\begin{vmatrix}
..&...&q_1&...&...\\ 
..&...&q_2&...&...\\ 
..&...&...&...&...\\ 
..&...&q_\ell&...&...\\ 
\end{vmatrix}=\frac{cQ(\mathcal{A})}{\alpha_{i_1}.\alpha_{i_2}...\alpha_{i_\ell}}.$$
This shows that the associated vector field $(q_1,\cdots,q_\ell)$ has no critical points in the complement.
\owari\\

\begin{cor}[Terao Theorem]
Let $\mathcal{A}$, $\mathcal{A}'$ be two arrangements with $L(\mathcal{A})\cong L(\mathcal{A}')$. If $\mathcal{A}$ is free then $\mathcal{A}'$ is also free. 
\end{cor}
\textbf{Proof}. 
We will prove that there are derivations $\theta_E,\theta_2',\cdots,\theta_i',\cdots,\theta_\ell' $ of $D(\mathcal{A}')$ such that they are linearly independent.\\
As in remark $3.4$, $D(\mathcal{A})$ has a basis $\{\theta_E,\theta_2,\cdots,\theta_i,\cdots,\theta_\ell \}$ such that the associated vector fields of every $\theta_i$ have no critical points in the complement. Thus, $\mathscr{C}_{\theta_i} (\mathcal{A})$ has no hidden constraints for every $i\geq 2$. It is followed from theorem $3.3$ that there are derivations $\theta_E,\theta_2',\cdots,\theta_i',\cdots,\theta_\ell' $ of $D(\mathcal{A}')$ having the same combinatorial structures.\\
By changing variables, we can assume that $v=(1,\cdots,1)\in M(\mathcal{A})$ is not a critical point of any associated vector fields of  derivations $\theta_i$ ( or $\theta_i' ), i\geq 2$. It follows that each of $\bar{\theta}_i\big|_v$ belongs to $(\bigcap_{j\in J_i}H_j)\setminus (\bigcup_{k\notin J_i} H_k).$\\
Since $\{\theta_E,\theta_2,\cdots,\theta_i,\cdots,\theta_\ell \}$ is a basis, $\{\bar{\theta}_E\big|_v,\bar{\theta}_2\big|_v,\cdots,\bar{\theta}_i\big|_v,\cdots,\bar{\theta}_\ell\big|_v \}$ is linearly independent. By choosing suitable polynomials $p$ as in remark $3.4$, we can assume that $(\bigcap_{j\in J_i}H_j)\setminus (\bigcup_{k\notin J_i} H_k)$ is different from $(\bigcap_{j\in J_t}H_j)\setminus (\bigcup_{k\notin J_t} H_k)$ for $i\neq t.$ Hence, $\{\bar{\theta}_E\big|_v,\bar{\theta}_2'\big|_v,\cdots,\bar{\theta}_\ell'\big|_v \}$ is linearly independent. \owari\\

\vspace{5mm}

\noindent
Tran Quoc Cong\\
Vietnamese Institute of Mathematics \\
18 Hoang Quoc Viet Road, Cau Giay District, Hanoi, Viet Nam\\
tqcong@gmail.com


\begin{thebibliography}{9}

\bibitem{ADS} T. Abe, A. Dimca and G. Sticlaru, Addition–deletion results for the minimal degree of logarithmic derivations of hyperplane arrangements and maximal Tjurina line arrangements. \textit { J. Algebraic Combin}, \textbf{54} (2020), 739-766.  

\bibitem{BC} M. Barakat and M. Cuntz,  Coxeter and crystallographic arrangements are inductively free. \textit {Adv. Math}, \textbf{229} (2012), 691-709.

\bibitem{Eis95} David Eisenbud,  \textit{Commutative algebra}, Graduate Texts in Mathematics, vol. \textbf{150}, Springer-
Verlag, New York, 1995, With a view toward algebraic geometry. MR 1322960 (97a:13001)
\bibitem{DS} G. Denham and A. Steiner, Geometry of logarithmic derivations of hyperplane arrangements. \textit { Contemporary Mathematics}, \textbf{790} (2023), 55-97. 

\bibitem{OT} P. Orlik and H. Terao, \textit{Arrangements of hyperplanes}.
Grundlehren der Mathematischen Wissenschaften, 
\textbf{300}. Springer-Verlag, Berlin, 1992.

\bibitem{Sa1} K. Saito, Theory of logarithmic differential forms and logarithmic vector fields.\textit{J. Fac. Sci. Univ. Tokyo, Sect. IA, Math.}
 \textbf{27}, 265–291 (1980)

\bibitem{T0}
H. Terao, Arrangements of hyperplanes and their freeness I, II. 
\textit{J. Fac. Sci. Univ. Tokyo}, \textbf{27} (1980), 293--320.   


 \bibitem{T2} H. Terao, Generalized exponents of a free arrangement of hyperplanes 
and Shepherd-Todd-Brieskorn formula. \textit{Invent. Math}, 
\textbf{63} (1981), 159--179. 

\bibitem{T3} H. Terao, Multiderivations of Coxeter arrangements. \textit{Invent. Math}. \textbf{148}, 659–674
(2002)
\bibitem{Yo}M. Yoshinaga, On the freeness of 3-arrangements. \textit{Bull. London Math. Soc.}
\textbf{37}, (2005), no. 1, 126-134.

\bibitem{Yul} S. Yuzvinsky, Cohomology of local sheaves on arrangement lattices. \textit{Proc.
Amer. Math. Soc.}\textbf{112} (1991), 1207-1217.

\bibitem{Yu2} S. Yuzvinsky, The first two obstructions to the freeness of arrangements.\textit{Trans. Amer. Math. Soc.} \textbf{335} (1993), 231-244.
\bibitem{Yu3}S. Yuzvinsky, Free and locally free arrangements with a given intersection
lattice. \textit{Proc. Amer. Math. Soc} \textbf{118} (1993), 745-752.
\bibitem{Zi} G. Ziegler, Combinatorial construction of logarithmic differential forms. \textit {Adv. in Math}, \textbf{76} (1989), 116-154.
\bibitem{Z0} G. Ziegler, Matroid Representations and free Arrangements. \textit{Trans. Amer. Math. Soc.}, \textbf{320} (1990), 525-541.

\bibitem{Z}G. M. Ziegler, Multiarrangements of hyperplanes and their freeness. in Singularities
(Iowa City, IA, 1986), 345-359, \textit { Contemp. Math,} \textbf{90}, Amer. Math.
Soc, Providence, RI, 1989.
\end{thebibliography}
\end{document}